\documentstyle[12pt]{article}
\title{The Hamilton-Waterloo problem for Hamilton cycles and $C_{4k}$-factors
\footnote{Research supported by National
Natural Science Foundation of China under Grant 10771137}} \vskip 0.5cm
\author{\small Hongchuan Lei\footnote{
\small{\it E-mail: chuan@sjtu.edu.cn}}\ \ \ \ \ \
\ \ \ \ Hung-Lin Fu\footnote{\small{\it E-mail:
hlfu@math.nctu.edu.tw}}\ \ \ \ \ \
\ \ \ \ Hao Shen\footnote{\small{\it Corresponding author. E-mail:
haoshen@sjtu.edu.cn}}\\
\small{\it  $^{\dag \S}$Department of Mathematics,\  Shanghai Jiao Tong University}\\
\small{\it  $^{\ddag}$Department of Applied Mathematics,\  National Chiao Tung University}\\
 } \topmargin -12mm \textheight
181mm \textwidth 114mm
\date{}
\begin{document}
\maketitle
\begin{minipage}{100mm}
{\small\bf Abstract}

\hbox{} \hskip 6mm In this paper we give a complete solution to the Hamilton-Waterloo
problem for the case of Hamilton cycles and $C_{4k}$-factors for all
positive integers $k$. \\

\end{minipage}\\

\topskip 1cm \textheight 8in

\section{\large\bf{Introduction}}
The Hamilton-Waterloo problem is a generalization of the
well known Oberwolfach problem, which asks for a 2-factorization of
the complete graph $K_n$ in which $r$ of its 2-factors are
isomorphic to a given 2-factor $R$ and s of its 2-factors are
isomorphic to a given 2-factor $S$ with $2(r+s)=n-1$. The most
interesting case of the Hamilton-Waterloo problem is that $R$ consists of cycles of length $m$ and
$S$ consists of cycles of length $k$, such a 2-factorization of
$K_n$ is called uniform and denoted by $HW(n;r,s;m,k)$. The
corresponding Hamilton-Waterloo problem is the problem for the
existence of an $HW(n;r,s;m,k)$.

There exists no 2-factorization of $K_n$ when $n$ is even since the degree of each vertex is odd. In this case, we consider the 2-factorizations of $K_n-I_n$(where $I_n$ is a 1-factor of $K_n$) instead. The corresponding 2-factorization is also denoted by $HW(n;r,s;m,k)$. Obviously $2(r+s)=n-2$.

It is easy to see that the following conditions are necessary for
the existence of an $HW(n;r,s;m,k)$:

{\bf Lemma 1.1.}
If there exists an $HW(n;r,s;m,k)$, then

\ \ $n\equiv 0\pmod{m}$ when $s=0$;

\ \ $n\equiv 0\pmod{k}$ when $r=0$;

\ \ $n\equiv 0\pmod{m}$ and $n\equiv 0\pmod{k}$ when $r\neq 0$ and
$s\neq 0$;

The Hamilton-Waterloo problem attracts much attention and progress
has been made by several authors. Adams, Billington, Bryant and
El-Zanati \cite{1} deal with the case $(m,k)\in
\{(3,5),(3,15),(5,15)\}$. Danziger, Quattrocchi and Stevens\cite{3}
give an almost complete solution for the case $(m,k)=(3,4)$, which
is stated below:

{\bf Theorem 1.2.} \cite{3} An $HW(n;r,s;3,4)$ exists if and only
if\\ $n\equiv 0 \pmod {12}$ and $(n,s)\neq (12,0)$ with the
following possible exceptions:

\ \ $n=24$ and $s=2,4,6$;

\ \ $n=48$ and $s=6,8,10,14,16,18$.

The case $(m,k)=(n,3)$, i.e. Hamilton cycles and triangle-factors,
is studied by Horak, Nedela and Rosa \cite{8}, Dinitz and
Ling \cite{4,5} and the following partial result
obtained:

{\bf Theorem 1.3.}
\cite{4,5,8}

(a) If $n\equiv 3 \pmod {18}$, then an $HW(n;r,s;n,3)$ exists except
possibly when $n=93,111,129,183,201$ and $r=1$;

(b) If $n\equiv 9 \pmod
{18}$, then an $HW(n;r,s;n,3)$ exists except $n=9$ and $r=1$, except possibly when $n=153,207$ and $r=1$;

(c) If $n\equiv 15 \pmod {18}$ and $r\in \{1,\frac{(n+3)}{6},\frac{(n+3)}{6}+2,\frac{(n+3)}{6}+3,\ldots,\frac{(n-1)}{2}\}$, then an $HW(n;r,s;n,3)$ exists except possibly when $n=123,141,159,177,\\213,249$ and
$r=1$.

For $n\equiv 0\pmod{6}$, the problem for the existence of an
$HW(n;r,\\s;n,3)$ is still open.

The cases $(m,k)\in \{(t,2t)|t>4\}$ and $(m,k)\in\{(4,2t)|t>3\}$
have been completely solved by Fu and Huang \cite{6}.

{\bf Theorem 1.4.}\cite{6}

(a) Suppose $t\geq 4$, an $HW(n;r,s;t,2t)$ exists if and only if $n\equiv
0 \pmod {2t}$.

(b) For an integer $t\geq 3$, an $HW(n;r,s;4,2t)$ exists if and only if $n\equiv
0 \pmod {4}$ and $n\equiv 0\pmod {2t}$.

For $r=0$ or $s=0$, the Hamilton-Waterloo problem is in fact the
problem for the existence of resolvable cycle decompositions of the
complete graph, which has been completely solved by Govzdjak
\cite{7}.

{\bf Theorem 1.5.}\cite{7} There exists a resolvable $m$-cycle
decomposition of $K_n$(or $K_n-I$ when n is even) if and only if
$n\equiv0\pmod{m}$, \\$(n,m)\neq(6,3)$ and $(n,m)\neq(12,3)$.

The purpose of this paper is to give a complete solution to the
Hamilton-Waterloo problem for the case of Hamilton cycles and
$C_{4k}$-factors which is stated in the following theorem.

{\bf Theorem 1.6.} For given positive integer $k$, an
$HW(n;r,s;n,4k)$ exists if and only if $r+s=[\frac{n-1}{2}]$ and
$n\equiv 0\pmod{4k}$ if $s>0$ or $n\geq 3$ if $s=0$.

\section{\large\bf  Preliminaries}
In this section, we provide some basic constructions.

For convenience, we introduce the following notations first. A
$C_m$-factor of $K_n$ is a spanning subgraph of $K_n$ in which each
component is a cycle of length $m$. Let $r+s= {{[(n - 1)}
\mathord{\left/
 {\vphantom {{[(n - 1)} 2}} \right.
 \kern-\nulldelimiterspace} 2}]
$ and
$$HW^*(n;m,k)=\{r|an\ HW(n;r,s;m,k)\ exists\}.$$
We use HC to
represent Hamilton cycle for short.

By Lemma 1.1, the necessary condition for the existence of
$HW(n;\\r,s;n,4k)$ with $s>0$ is $n\equiv 0\pmod{4k}$, we assume
$n=4kt$ and the vertex set of $K_n$ is $Z_{2t}\times Z_{2k}$. We
write $V_i=\{i\}\times Z_{2k}=\{i_0,i_1,\ldots,i_{2k-1}\}$ for $i\in
Z_{2t}$. Let $K_{V_i,V_j}$ be the complete bipartite graph define on two partite sets $V_i$ and $V_j$, and $K_{V_i}$ be the complete graph of
order $2k$ define on the vertex set $V_i$. Obviously,
\[E(K_{4kt})=\bigcup\limits_{i=0}^{2t-1}{E(K_{V_i})}\cup
\bigcup\limits_{i\neq j}^{}E(K_{V_i,V_j}).\]

Further for $d\in Z_{2k}$, we define sets of edges $(i,j)_d=\{(i_lj_{l+d})|l\in Z_{2k}\}$ for $i,j\in Z_{2t}$. Clearly, $(i,j)_d$ is a perfect matching in $K_{V_i,V_j}$. In fact, $$E(K_{V_i,V_j})=\bigcup\limits_{d=0}^{2k-1}{(i,j)_d}.$$

The following lemmas are useful in our constructions.

{\bf Lemma 2.1.}
\cite{6} Let $I_{2n}=\{(v_0v_n)\}\cup\{(v_iv_{2n-i})|1\leq
i\leq n-1\}$. Then $K_{2n}-I_{2n}$ can be
decomposed into $n-1$ HCs, Each HC can be decomposed into two 1-factors. Moreover, by reordering the vertices of $K_{2n}$ if necessary, we may assume one of the HCs is $(v_0,v_1,\ldots,
v_{2n-1})$.

The following lemma is a generalization of Lemma 1 in \cite{8}.

{\bf Lemma 2.2.} Let $\pi$ be a permutation of $Z_{2t}$,
$d_0,d_1,\ldots,d_{2t-1}$ be nonnegative integers. Then the set of
edges \[(\pi(0),\pi(1))_{d_0}\cup(\pi(1),\pi(2))_{d_1}\cup\cdots
\cup (\pi(2t-1),\pi(0))_{d_{2t-1}}\] forms an HC of $K_n$ if
$d_0+d_1+\cdots +d_{2t-1}$ and $2k$ are relatively prime.

{\bf Proof.}
Set $d=d_0+d_1+\cdots +d_{2t-1}$, then arrange the edges as
\[H=(\pi(0)_0,\pi(1)_{d_0},\pi(2)_{d_0+d_1},\cdots,
\pi(0)_{d},\pi(1)_{d+d_0},\cdots,\pi(2t-1)_{2kd-d_{2t-1}}).\]
Since $(d,2k)=1$, the vertices
\[\pi(i)_{d_0+d_1+\cdots+d_{i-1}},\pi(i)_{d+d_0+d_1+\cdots+d_{i-1}},\ldots,\pi(i)_{(2k-1)d+d_0+d_1+\cdots+d_{i-1}}\] are mutually distinct for $i\in Z_{2t}$. Thus all vertices in $H$ are mutually distinct, so $H$ is an HC.
$\Box$

{\bf Lemma 2.3.}
Let $d_1,d_2$ be nonnegative integers. If $d_1-d_2$ and $2k$ are
relatively prime, then the set of edges $(i,j)_{d_1}\cup (i,j)_{d_2}$
forms a cycle of length $4k$ on the vertex set $V_i\cup V_j$.

{\bf Proof.}
It's a direct consequence of Lemma 2.2. Arranging the edges as a cycle
$(i_0,j_{d_1},i_{d_1-d_2},j_{2d_1-d_2},\cdots,j_{2kd_1-(2k-1)d_2})$
completes the proof.$\Box$

\section{\large\bf{Proof of the main theorem}}
With the above preparations, now we are ready to prove our main theorem.

Let $\widetilde{G}$ be a complete graph defined on $\{V_0,V_1,\ldots,V_{2t-1}\}$. By Lemma 2.1, $\widetilde{G}$ can be decomposed into $2t-1$
1-factors, denoted by $\widetilde{F}_1,\widetilde{F}_2,\ldots,\widetilde{F}_{2t-1}$, and $\widetilde{F}_{2i-1}\cup
\widetilde{F}_{2i}$ forms an HC for $i=1,2,\ldots,t-1$. By reordering the
vertices if necessary, we may assume
\[\widetilde{F}_1=\{V_0V_1,V_2,V_3,\ldots,V_{2t-2}V_{2t-1}\},\]
\[\widetilde{F}_2=\{V_1V_2,V_3V_4,\ldots,V_{2t-1}V_0\},\]
\[\widetilde{F}_{2t-1}=\{V_0V_t\}\cup\{V_iV_{2t-i}|i=1,2,\ldots,t-1\}.\]
Let
\[F_x=\bigcup\limits_{{V_i}{V_j} \in E({{\widetilde F}_x})} {E({K_{{V_i},{V_j}}})}\ \ for\ x\in Z_{2t}\backslash\{0\}\]
and
\[H_l=(0,1)_l\cup (1,2)_{2k-l}\cup (2,3)_l\cup \cdots\cup (2t-1,0)_{2k-l} \ \ for \ l\in Z_{2k}.\]
Then $F_1\cup F_2=H_0\cup
H_1\cup\cdots\cup H_{2k-1}$.

{\bf Lemma 3.1.}
$F_{2i-1}\cup F_{2i}(i=0,1,\ldots,k-1)$ can be decomposed into
$r_i\in \{0,2,\ldots,2k\}$ HCs and $2k-r_i$ $C_{4k}$-factors of
$K_{n}$.

{\bf Proof.}
We only give the proof for the case $i=1$, i.e. $F_1\cup F_2$, the
remaining cases are similar.

For $l=0,1,\ldots,k-1$, $H_{2l}\cup H_{2l+1}$ can be decomposed into
two edge sets:
\[\bigcup\limits_{j = 0}^{t - 1} {( {{(2j,2j + 1)}_{2l}}\bigcup {{{(2j,2j + 1)}_{2l + 1}}} ) ,} \]
\[\bigcup\limits_{j = 0}^{t - 1} {( {{(2j + 1,2j + 2)}_{2k - 2l}}\bigcup {{{(2j + 1,2j + 2)}_{2k - 2l - 1}}} ) }, \]
by Lemma 2.3, each forms a $C_{4k}$-factor of $K_n$.

Similarly, $H_{2l}\cup H_{2l+1}$ can be decomposed into another two
edge sets:
\[(H_{2l}-(2t-1,0)_{2k-2l})\cup (2t-1,0)_{2k-2l-1},\]
\[(H_{2l+1}-(2t-1,0)_{2k-2l-1})\cup (2t-1,0)_{2k-2l},\]
 by Lemma 2.2,
each forms an HC of $K_{n}$.

Finally, by decomposing $H_{2l}\cup H_{2l+1}$
into two HCs when $l\in\{0,1,\ldots,\frac{r_i}{2}-1\}$ or into two $C_{4k}$-factors when $l\in\{\frac{r_i}{2},\frac{r_i}{2}+1,\ldots,k-1\}$, we have the proof.$\Box$

{\bf Lemma 3.2.}
For each $i\in Z_{2t}\backslash\{0\}$, $F_{i}\cup(\bigcup\limits_{i\in Z_{2t}}^{}{K_{V_i}})$ can be decomposed into $2k-1$ $C_{4k}$-factors and a 1-factor of $K_{n}$.

{\bf Proof.}
Noticing that
$F_{i}\cup(\bigcup\limits_{i\in Z_{2t}}^{}{K_{V_i}})=tK_{4k}$ and these
complete graphs of order $4k$ are edge-disjoint. By Lemma 2.1,
each can be decomposed into $2k-1$ HCs and one 1-factor of $K_{4k}$.
Hence, these HCs and 1-factors form $2k-1$ $C_{4k}$-factors and a
1-factor of $K_{n}$. This concludes the proof.
$\Box$

For convenience in presentation, we use ${\rm X}$ to denote $\bigcup\limits_{i\in Z_{2t}}{K_{V_i}}$ in what follows.

{\bf Proposition 3.3.}
$\{0,2,4,\ldots,\frac{n}{2}-2k\}\subseteq
HW^*(n;n,4k)$ for all positive integers $n\equiv 0 \pmod{4k}$.

{\bf Proof.}
Since $K_{n}=F_1\cup F_2\cup \cdots \cup
F_{2t-1}\cup {\rm X}$, applying Lemma
3.2 to $F_{2t-1}\cup{\rm X}$ and
Lemma 3.1 to $F_{2i}\cup F_{2i-1}(1\leq i\leq t-1)$ completes the
proof.
$\Box$

{\bf Proposition 3.4.}
$\{1,3,5,\ldots,\frac{n}{2}-4k+1\}\subseteq HW^*(n;n,4k)$ for all
positive integers $n\equiv 0 \pmod{4k}$.

{\bf Proof.}
First, by Lemma 3.2, we decompose
$F_{2}\cup{\rm X}$ into $2k-1$
$C_{4k}$-factors and a 1-factor. Without loss of generality,
assume the 1-factor is $I_n^{'}=(1,2)_0\cup (3,4)_0\cup \cdots\cup
(2t-1,0)_0$.

Since $E({F_1}) = \bigcup\limits_{i = 0}^{2k - 1}
{({{(0,1)}_i}\bigcup {{{(2,3)}_i}} }  \cdots {(2t - 2,2t - 1)_i})$,
we decompose $E({F_1})\cup I_n^{'}$ into $k-1$ $C_{4k}$-factors, an HC
and a 1-factor:
\[C_i=((0,1)_{2i-1}\cup(0,1)_{2i})\cup((2,3)_{2i-1}\cup(2,3)_{2i})\cup\cdots\cup ((2t-2,2t-1)_{2i-1}\cup\]
\[(2t-2,2t-1)_{2i}),\ \
i=1,2,\ldots,k-1,\]
\[HC_1=(0,1)_{2k-1}\cup(1,2)_0\cup(2,3)_0\cup\cdots\cup(2t-2,2t-1)_0,\]
\[I_n=(0,1)_0\cup (2,3)_{2k-1}\cup (4,5)_{2k-1}\cdots\cup (2t-2,2t-1)_{2k-1}.\]
It is straightforward to verify that $C_i$ is a $C_{4k}$-factor,
$HC_1$ is an HC, $I_n$ is a 1-factor and they are edge-disjoint.

Finally, applying Lemma 3.1 to $F_{2i-1}\cup F_{2i}(2\leq i\leq
t-1)$ gives $\{1,3,5,\ldots,\frac{n}{2}-4k+1\}\subseteq
HW^*(n;n,4k)$.
$\Box$

{\bf Lemma 3.5.}
If $r_1\in\{2k,2k+1,2k+2,\ldots,4k-1\}$, then $F_{1}\cup F_{2}\cup
F_{2t-1}\cup{\rm X}$ can be
decomposed into $r_1$ HCs, $4k-1-r_1$ $C_{4k}$-factors and a
1-factor of $K_{n}$.

{\bf Proof.}
It is well known that every complete graph with even order can be
decomposed into Hamilton paths\cite{2}. Noticing that
\[F_{2t-1}\cup{\rm X}=\{K_{V_0\cup
V_t}\}\cup \{K_{V_i\cup V_{2t-i}}|i=1,2,\ldots,t-1\}=tK_{4k}\] and these
complete graphs of order $4k$ have no common vertex. Let
$P_{i,j}[u\ldots v]$ be the Hamilton path of $K_{V_{i}\cup V_{j}}$
with $u$ and $v$ as its end vertices. We may decompose
$F_{2t-1}\cup{\rm X}$ into
$\{P_0,P_1,\ldots,P_{2k-1}\}$ where
\[P_j=\{P_{0,t}[0_j,\ldots,t_j]\}\cup \{P_{i,2t-i}[i_j,\ldots,(2t-i)_j]|i=1,2,\ldots,t-1\}.\]
For each $j$, connecting the Hamilton paths of $P_j$ with $t$ edges $(0_j1_j),\\ (2_j3_j),\ldots, ((2t-2)_j(2t-1)_j)\in (0,1)_0\cup(2,3)_0\cup\cdots\cup(2t-2,2t-1)_0\subseteq H_0$ which gives an $HC$. Then we have $2k$ Hamilton cycles $HC_j$, $j\in Z_{2k}$, when $t$ is odd,
\[\begin{array}{*{20}{c}}
   {H{C_j} = } \hfill & {({0_j},{1_j},{P_{1,2t - 1}}[{1_j}, \ldots ,{{(2t - 1)}_j}],{{(2t - 1)}_j},{{(2t - 2)}_j},} \hfill  \\
   {} \hfill & {{P_{2t - 2,2}}[{{(2t - 2)}_j}, \ldots ,{2_j}], \ldots ,{{(t - 1)}_j},{t_j},{P_{t,0}}[{t_j}, \ldots ,{0_j}]);} \hfill  \\
\end{array}\]
when $t$ is even,
\[\begin{array}{*{20}{c}}
   {H{C_j} = } \hfill & {({0_j},{1_j},{P_{1,2t - 1}}[{1_j}, \ldots ,{{(2t - 1)}_j}],{{(2t - 1)}_j},{{(2t - 2)}_j},} \hfill  \\
   {} \hfill & {{P_{2t - 2,2}}[{{(2t - 2)}_j}, \ldots ,{2_j}], \ldots ,{{(t + 1)}_j},{t_j},{P_{t,0}}[{t_j}, \ldots ,{0_j}]).} \hfill  \\
\end{array}\]

Then we can decompose $H_1\cup
(H_0-(0,1)_0\cup(2,3)_0\cup\cdots\cup(2t-2,2t-1)_0)$ into an HC and
a 1-factor, or a $C_{4k}$-factor and a 1-factor. In the first
case, let
\[HC_{2k}=H_1\cup(2t-1,0)_{0}-(2t-1,0)_{2k-1},\]
\[I_n=(1,2)_0\cup (3,4)_0\cup\cdots\cup(2t-3,2t-2)_0\cup(2t-1,0)_{2k-1}.\]By Lemma
2.2, $HC_{2k}$ forms an HC. $I_n$ is a 1-factor. In the second case,
let
\[C=\bigcup\limits_{j = 0}^{t - 1} {\{ {{(2j + 1,2j + 2)}_0}\bigcup {{{(2j + 1,2j + 2)}_{2k - 1}}} \} },\]
\[I_n^{'}=(0,1)_1\cup (2,3)_1\cup\cdots\cup(2t-2,2t-1)_{1}.\]
By Lemma 2.3, $C$ is a $C_{4k}$-factor and $I_n^{'}$ is a 1-factor.

Finally, in the same way as Lemma 3.1, for each $r_1\in\{2k,2k+2,2k+4,\ldots,4k-2\}$, we decompose each $H_{2l}\cup
H_{2l+1}$ into two HCs for $l\in\{1,2,\ldots,\frac{r_1}{2}\}$ or two $C_{4k}$-factors for $l\in\{\frac{r_1}{2}+1,\frac{r_1}{2}+2,\ldots,k-1\}$. Then we have the proof.
$\Box$

{\bf Proposition 3.6.}
$\{2k,2k+1,2k+2,\ldots,\frac{n-2}{2}\}\subseteq HW^*(n;n,4k)$ for
all positive integers $n\equiv 0 \pmod{4k}$.

{\bf Proof.}
Let $r=p\cdot 2k+q,$ where $0\leq q<2k$. If $2k\leq r\leq 2kt-2k$ and $q$ is even, by Lemma 3.5, we may decompose $F_{1}\cup F_{2}\cup F_{2t-1}\cup{\rm X}$ into $2k$ HCs, $2k-1$ $C_{4k}$-factors and a 1-factor. By Lemma 3.1, we may decompose $F_{2i-1}\cup F_{2i}$ into $2k$ HCs for each $2\leq i\leq p$, $F_{2p+1}\cup F_{2p+2}$ into $q$ HCs and $2k-q$ $C_{4k}$-factors, and $F_{2j-1}\cup F_{2j}$ into $2k$ $C_{4k}$-factors for each $p+2\leq j\leq t-1$. Then we have \[\{2k,2k+2,\ldots,2kt-2k\}\subseteq HW^*(n;n,4k).\]

If $2k\leq r\leq 2kt-2k$ and $q$ is odd, by Lemma 3.5, we may decompose $F_{1}\cup F_{2}\cup
F_{2t-1}\cup{\rm X}$ into $2k+1$ HCs, $2k-2$ $C_{4k}$-factors and a 1-factor. By Lemma 3.1, we may decompose $F_{2i-1}\cup F_{2i}$ into $2k$ HCs for each $2\leq i\leq p$, $F_{2p+1}\cup F_{2p+2}$ into $q-1$ HCs and $2k-q+1$ $C_{4k}$-factors, and $F_{2j-1}\cup F_{2j}$ into $2k$ $C_{4k}$-factors for each $p+2\leq j\leq t-1$. Then we have \[\{2k+1,2k+3,\ldots,2kt-2k-1\}\in HW^*(n;n,4k).\]

If $2kt-2k<r\leq \frac{n-2}{2}$ and $q$ is even, by Lemma 3.5, we may decompose $F_{1}\cup F_{2}\cup F_{2t-1}\cup{\rm X}$ into $4k-2$ HCs, a $C_{4k}$-factor and a 1-factor. When $q+2<2k$, by Lemma 3.1, we may decompose $F_{2i-1}\cup F_{2i}$ into $2k$ HCs for each $2\leq i\leq p-1$, $F_{2p-1}\cup F_{2p}$ into $q+2$ HCs and $2k-q-2$ $C_{4k}$-factors, and $F_{2j-1}\cup F_{2j}$ into $2k$ $C_{4k}$-factors for each $p+1\leq j\leq t-1$; when $q+2=2k$, we decompose $F_{2i-1}\cup F_{2i}$ into $2k$ HCs for each $2\leq i\leq p$ and $F_{2j-1}\cup F_{2j}$ into $2k$ $C_{4k}$-factors for each $p+1\leq j\leq t-1$. Then we have \[\{2kt-2k+2,2kt-2k+4,\ldots,2kt-2\}\in HW^*(n;n,4k).\]

If $2kt-2k<r\leq \frac{n-2}{2}$ and $q$ is odd, by Lemma 3.5, we may decompose $F_{1}\cup F_{2}\cup F_{2t-1}\cup{\rm X}$ into $4k-1$ HCs and a 1-factor. When $q+1=2k$,by Lemma 3.1, we may decompose each $F_{2i-1}\cup F_{2i}$into $2k$ HCs for each $2\leq i\leq p$ and $F_{2j-1}\cup F_{2j}$ into $2k$ $C_{4k}$-factors for each $p+1\leq i\leq t-1$; when $q+1\neq 2k$,  we decompose $F_{2i-1}\cup F_{2i}$into $2k$ HCs for each $2\leq i\leq p-1$, $F_{2p-1}\cup F_{2p}$ into $q+1$ HCs and $2k-q-1$ $C_{4k}$-factors, and $F_{2j-1}\cup F_{2j}$ into $2k$ $C_{4k}$-factors for each $p+1\leq j\leq t-1$. Then we have \[\{2kt-2k+1,2kt-2k+3,\ldots,2kt-1\}\in HW^*(n;n,4k).\Box\]

Combining Proposition 3.3, Proposition 3.4 and Proposition 3.6, we have the main
result of this paper.

{\bf Theorem 3.7.}
$\{0,1,2,\ldots,\frac{n-2}{2}\}=HW^*(n;n,4k)$ for all positive
integers $n\equiv 0 \pmod{4k}$.

{\bf Proof.} For $n=4k$, the theorem is obvious by Theorem 1.5. For
$n=8k$, the result is also correct by Theorem 1.4. When $n>8k$, we
have $\frac{n}{2}-2k>2k$ and $\frac{n}{2}-4k+1\geq 2k+1$, then
combining with Proposition 3.3, Proposition 3.4 and Proposition 3.6
completes the proof. $\Box$

\section{\large\bf{Concluding remarks}}

It would be interesting to determine the necessary and sufficient
conditions for the existence of an $HW(n;r,s;n,k)$ for any even
integer $k$. As a first step, we proved in this paper that for any
integer $k\equiv 0\pmod 4$ the necessary condition for the existence
of $HW(n;r,s;n,k)$ is $n\equiv 0\pmod{k}$, and the necessary
condition is also sufficient. The next step is for the case when
$k\equiv 2\pmod 4$, we conjecture that for $k\equiv 2\pmod 4$ and
$s>0$ there exists an $HW(n;r,s;n,k)$ if and only if $n\equiv 0\pmod
k$.


\begin{thebibliography}{1}
\bibitem{1}P. Adams, E.J. Billington, D.E. Bryant, S.I.
El-Zanati, On the Hamilton-Waterloo problem, Graph Combin. 18 (2002)
31-51.
\bibitem{2}C. J. Colbourn, J. H. Dinitz (Editors), The CRC Handbook of Combinatorial Designs.
2nd edn, CRC Press Series on Discrete Mathematics, CRC, Boca Raton,
2007.
\bibitem{3}P. Danziger, G. Quattrocchi, B. Stevens, The Hamilton-Waterloo problem for
cycle sizes 3 and 4, J. Combin. Designs. 17 (2009) 342-352.
\bibitem{4}J. H. Dinitz, A. C. H. Ling, The Hamilton-Waterloo
problem with triangle-factors and Hamilton cycles: The case $n\equiv
3\pmod {18}$, J. Combin. Math. Combin. Comput., to appear.
\bibitem{5}J.H. Dinitz, A. C. H. Ling, The Hamilton-Waterloo
problem: the case of triangle-factors and one Hamilton cycle, J.
Combin. Designs. 17 (2009) 160-176.
\bibitem{6}H. L. Fu, K. C. Huang, The Hamilton-Waterloo problem for
two even cycles factors, Taiwanese Journal of Mathematics 12 (2008)
933-940.
\bibitem{7}P. Govzdjak, On the Oberwolfach problem for the
complete mutigraphs, Discrete Math. 173(1997)61-69.
\bibitem{8}P. Horak, R. Nedela, A. Rosa, The Hamilton-Waterloo
problem: the case of Hamilton cycles and triangle-factors, Discrete
Math. 284 (2004) 181-188.
\end{thebibliography}
\end{document}